\documentclass[12pt,a4paper]{article}
\usepackage{latexsym,amssymb,amsmath,mathrsfs}

\usepackage{sectsty}
\usepackage{color}
\usepackage{bm}
\usepackage[T1]{fontenc}
\usepackage[anchorcolor=blue,%
 bookmarks=true,%
 bookmarksnumbered=true,%
]{hyperref}
\usepackage{cite}
\usepackage{pst-all}
\usepackage{graphicx}

\allsectionsfont{\normalsize\sc\center}

\newtheorem{thm}{Theorem}[section]
\newtheorem{prop}[thm]{Proposition}

\newtheorem{lem}[thm]{Lemma}

\newtheorem{remark}[thm]{Remark}

\newenvironment{pf}{\par\begin{trivlist}%
\item[]{\bf Proof.}\ }{\hfill $\square$ \end{trivlist}\par}

\makeatletter \@addtoreset{equation}{section} \makeatother

\newcommand{\R}{\mathbb{R}}

\DeclareMathOperator{\Ric}{Ric}

\newcommand{\IR}{\operatorname{InRad}}
\newcommand{\Div}{\operatorname{div}}

\renewcommand{\d}{\mathrm{d}}

\newcommand{\s}{\mathfrak{s}}

\newcommand{\eps}{{\varepsilon}}

\newcommand{\inte}{\mathrm{Int}\,}

\title{\large\bf  Comparison geometry of manifolds with boundary under lower \boldmath$N$-weighted Ricci curvature bounds with 
\boldmath$\varepsilon$-range}
\author{Kazuhiro Kuwae\thanks{Department of Applied Mathematics, Fukuoka University,
Fukuoka 814-0180, Japan ({\sf kuwae@}  {\sf fukuoka-u.ac.jp}) . Supported in part by JSPS Grant-in-Aid for Scientific Research (KAKENHI) 17H02846 and by fund (No.:185001) from the Central Research Institute of Fukuoka University.}
\ \ and\ \ 
Yohei Sakurai\thanks{Advanced Institute for Materials Research, Tohoku University
2-1-1 Katahira, Aoba-ku, Sendai, 
980-8577, Japan ({\sf yohei.sakurai.e2@tohoku.ac.jp}). Supported in part by JSPS Grant-in-Aid for Scientific Research on Innovative Areas ``Discrete Geometric Analysis for Materials Design" 17H06460).
}
}
\date{}
\begin{document}
\maketitle
\begin{abstract}
We study comparison geometry of manifolds with boundary under a lower $N$-weighted Ricci curvature bound for $N\in ]-\infty,1]\cup [n,+\infty]$ with $\varepsilon$-range introduced by Lu-Minguzzi-Ohta \cite{LMO:CompaFinsler}.
We will conclude splitting theorems,
and also comparison geometric results for inscribed radius, volume around the boundary, and smallest Dirichlet eigenvalue of the weighted $p$-Laplacian.
Our results interpolate those for $N\in [n,+\infty[$ and $\varepsilon=1$, and for $N\in ]-\infty,1]$ and $\varepsilon=0$ by the second named author.
\end{abstract}


 
\section{Introduction}

\subsection{Background}
We first recall the notion of the weighted Ricci curvature,
and comparison geometry of manifolds without boundary under lower $N$-weighted Ricci curvature bounds.

Let $(M,g,f)$ be an $n$-dimensional weighted Riemannian manifold, namely, $(M,g)$ is an $n$-dimensional complete Riemannian
manifold, and $f\in C^{\infty}(M)$.
For $N\in]-\infty,\,+\infty]$, the {\it $N$-weighted Ricci curvature} is defined as follows (\cite{Ba}, \cite{BE1}, \cite{Lich}): 
\begin{align*}
\Ric_f^N:=\Ric_g+{\rm \nabla^2}f-\frac{df\otimes df}{N-n}.
\end{align*}
Here if $N=+\infty$,
the last term should be interpreted as the limit $0$,
and when $N=n$,
we only consider a constant function $f$ such that $\Ric_f^n:=\Ric_g$.

In the rest of this subsection,
the boundary $\partial M$ of $M$ is assumed to be empty.
In such a case,
under a classical curvature condition
\begin{align}
\Ric_f^N\geq Kg \label{eq:RicciLowerBdd}
\end{align}
for $K\in \mathbb{R}$ and $N\in [n,+\infty[$,
comparison geometry has been developed by \cite{FLZ}, \cite{Lo}, \cite{Qi}, \cite{WW},
and so on (see also \cite{FLL}, \cite{LL1}, \cite{LL2}, \cite{Li1}, \cite{Li2}).
In recent years,
comparison geometry for the complementary case of $N\in ]-\infty,n[$ has begun to be investigated (see e.g., \cite{KL}, \cite{KSa}, \cite{KS}, \cite{Lim}, \cite{LMO:CompaFinsler}, \cite{Mai1}, \cite{Mai2}, \cite{Mineg}, \cite{Oh<0}, \cite{Sak}, \cite{Wy}, \cite{WyYero}).
Wylie-Yeroshkin \cite{WyYero} have introduced a curvature condition
\begin{align}
\Ric_f^1\geq (n-1)\kappa e^{-\frac{4f}{n-1}}g \label{eq:RicciLowerBddWY}
\end{align}
for $\kappa \in \mathbb{R}$ in view of the study of weighted affine connection,
and obtained a Laplacian comparison theorem, Bonnet-Myers type diameter comparison theorem, Bishop-Gromov type volume comparison theorem, and rigidity results for the equality cases.
The first named author and Li \cite{KL} have provided a generalized condition
\begin{align}
\Ric_f^N\geq (n-N)\kappa e^{-\frac{4f}{n-N}}g\label{eq:RicciLowerBddKL}
\end{align}
with $N\in ]-\infty,1]$,
and extended the comparison geometric results in \cite{WyYero}.

Very recently,
Lu-Minguzzi-Ohta \cite{LMO:CompaFinsler} have introduced a new curvature condition that interpolates the conditions \eqref{eq:RicciLowerBdd} with $K=(N-1)\kappa$, \eqref{eq:RicciLowerBddWY} and \eqref{eq:RicciLowerBddKL}.
For $N\in ]-\infty,1]\cup [n,+\infty]$,
the notion of the \emph{$\varepsilon$-range} has played a crucial role in their work:
\begin{align}
\varepsilon=0\ \text{ for }\ N=1,\qquad \varepsilon\in ]-\sqrt{\varepsilon_0},\sqrt{\varepsilon_0}[  \ \text{ for }\  N\ne1,n,\qquad \varepsilon\in\R \ \text{ for } \ N=n,\label{eq:epsilonrange}
\end{align}
where we set
\begin{align*}
\varepsilon_0:=\frac{N-1}{N-n},
\end{align*}
and interpret it as the limit $1$ when $N=+\infty$.
In this $\eps$-range,
they have studied
\begin{align}
\Ric_f^N\geq c^{-1}\,\kappa\, e^{-\frac{4(1-\varepsilon)f}{n-1}}g,\label{eq:LMOLowerBddd}
\end{align}
where $c=c_{N,\varepsilon}\in]0,1]$ is a positive constant defined by 
\begin{align}
c:=\frac{1}{n-1}\left(1-\varepsilon^2 \frac{N-n}{N-1} \right)\label{eq:constant}
\end{align}
if $N\ne1$,
and $c:=(n-1)^{-1}$ if $N=1$. 
When $N\in [n,+\infty[$ and $\varepsilon=1$ with $c=(N-1)^{-1}$, the curvature condition \eqref{eq:LMOLowerBddd} is reduced to \eqref{eq:RicciLowerBdd} with $K=(N-1)\kappa$.
Further,
when $N=1$ and $\varepsilon =0$ with $c=(n-1)^{-1}$, it covers \eqref{eq:RicciLowerBddWY},
and when $N\in ]-\infty,1]$ and $\varepsilon=\varepsilon_0$ with $c=(n-N)^{-1}$, 
it does \eqref{eq:RicciLowerBddKL}. 
Under the condition \eqref{eq:LMOLowerBddd},
they have shown a Laplacian comparison theorem,
and concluded a diameter bound of Bonnet-Myers type,
and a volume bound of Bishop-Gromov type under density bounds.
The authors \cite{KSa} have examined rigidity phenomena for the equality cases in \cite{LMO:CompaFinsler}.

\subsection{Setting}

In this article,
we focus on the case where the boundary $\partial M$ is non-empty.
Under a lower Ricci curvature bound,
Heintze-Karcher \cite{HK}, Kasue \cite{K2}, \cite{K3} has studied comparison geometry of (unweighted) manifolds with boundary assuming a lower mean curvature bound for the boundary.
The second named author \cite{Sa2}, \cite{Sak} has developed that of weighted manifolds with boundary extending the work of \cite{HK}, \cite{K2}, \cite{K3}.

On a weighted Riemannian manifold with boundary $(M,g,f)$,
the \textit{weighted mean curvature} at $z\in \partial M$ is defined by
\begin{equation*}
H_{f,z}:=H_{z}+g(\nabla f,u_z),
\end{equation*}
where $H_z$ denotes the mean curvature at $z$ defined as the trace of the shape operator with respect to the unit inner normal vector $u_z$.
For a function $\Lambda:\partial M\to \mathbb{R}$,
we write $H_{f,\partial M}\geq \Lambda$ when $H_{f,z}\geq \Lambda(z)$ for all $z\in \partial M$.
The second named author \cite{Sa2} has developed comparison geometry under a curvature condition
\begin{equation}
\Ric_f^N\geq (N-1)\kappa g,\quad H_{f,\partial M}\geq (N-1)\lambda \label{eq:Sa2setting}
\end{equation}
for $\kappa,\lambda \in \mathbb{R}$ and $N\in [n,+\infty[$
that corresponds to \eqref{eq:RicciLowerBdd}.
After that,
inspired by the work of Wylie-Yeroshkin \cite{WyYero},
the second named author \cite{Sak} has also considered a condition
\begin{align}
\Ric_f^N\geq (n-1)\,\kappa\, e^{-\frac{4f}{n-1}}g,\quad H_{f,\partial M}\geq (n-1)\lambda e^{-\frac{2f}{n-1}} \label{eq:Saksetting}
\end{align}
for $N\in ]-\infty,1]$, and derived comparison geometric results.

The aim of this note is to establish comparison geometry under the curvature condition \eqref{eq:LMOLowerBddd},
which enables us to treat the results in \cite{Sa2}, \cite{Sak} in a unified way.
Our setting is as follows:
For $N\in ]-\infty,1]\cup [n,+\infty]$ and $\eps \in \mathbb{R}$ in the range \eqref{eq:epsilonrange},
\begin{align}
\Ric_f^N\geq c^{-1}\,\kappa\, e^{-\frac{4(1-\varepsilon)f}{n-1}}g,\quad H_{f,\partial M}\geq c^{-1}\lambda e^{-\frac{2(1-\varepsilon)f}{n-1}}, \label{eq:setting}
\end{align}
where $c$ is defined as \eqref{eq:constant}.
When $N\in [n,+\infty[$ and $\varepsilon=1$, the condition \eqref{eq:setting} is reduced to \eqref{eq:Sa2setting},
and when $N\in ]-\infty,1]$ and $\varepsilon =0$, it is reduced to \eqref{eq:Saksetting}.
Under \eqref{eq:setting},
we prove several splitting theorem (see Section \ref{sec:Splitting theorems}),
and also show comparison geometric results for inscribed radius (see Section \ref{sec:Inscribed radii}), volume of metric neighborhoods of the boundary (see Section \ref{sec:Volume growths}), and Dirichlet eigenvalue of the weighted $p$-Laplacian (see Section \ref{sec:Spectrum}).


\section{Preliminaries}\label{sec:Preliminaries}
In what follows,
let $(M,g,f)$ be an $n$-dimensional weighted Riemannian manifold,
and we always fix $N\in ]-\infty,1]\cup [n,+\infty]$ and $\eps \in \mathbb{R}$ in the range \eqref{eq:epsilonrange}.
We basically use the same notation and terminology as in \cite{Sak}.
We also assume familiarity with basic concepts in geometry of manifolds with boundary such as the cut locus for the boundary (cf. \cite{Sa1}, and Section 2 in \cite{Sa2}, \cite{Sak}).

In this section,
as a preliminary,
we recall a Laplacian comparison theorem for the distance function from a single point under the curvature condition \eqref{eq:LMOLowerBddd} obtained by Lu-Minguzzi-Ohta \cite{LMO:CompaFinsler},
and a rigidity results for the equality case investigated in \cite{KSa}.

\subsection{Laplacian comparison theorem from a single point}
Let $\inte M$ stand for the interior of $M$.
For $x \in \inte M$,
let $U_{x}M$ be the unit tangent sphere at $x$,
which can be identified with the $(n-1)$-dimensional unit sphere $(\mathbb{S}^{n-1},g_{\mathbb{S}^{n-1}})$.
Let $\rho_{x}:M\to \mathbb{R}$ be the distance function from $x$ defined as $\rho_{x}(y):=d(x,y)$,
where $d$ denotes the Riemannian distance function on $M$.
For $v \in U_{x} M$,
let $\gamma_{v}:[0,T]\to M$ be the unit speed geodesic with $\gamma_{v}(0)=x$ and ${\dot\gamma}_{v}(0)=v$. 
We define $\tau:U_{x}M\to ]0,+\infty]$ by
\begin{equation*}\label{eq:pointed cut value}
\tau(v):=\sup \left\{\, t>0 \mid \rho_{x}(\gamma_{v}(t))=t,\, \gamma_{v}([0,t[)\subset \inte M  \,\right\}.
\end{equation*}
Let $s_{f,v}:[0,\tau(v)]\to [0,+\infty]$ denote a function defined by
\begin{equation*}
s_{f,v}(t):=\int^{t}_{0}\,  e^{-\frac{2(1-\eps)f(\gamma_{v}(\xi))}{n-1}}\,\d \xi.
\end{equation*}

For $\kappa \in \mathbb{R}$,
we denote by $\mathfrak{s}_{\kappa}(s)$ the solution of the equation $\psi''(s)+\kappa \psi(s)=0$ with $\psi(0)=0$ and $\psi'(0)=1$.
We define
\begin{equation*}\label{eq:pointed model mean curvature}
H_{\kappa}(s):=-c^{-1}\frac{\mathfrak{s}'_{\kappa}(s)}{\mathfrak{s}_{\kappa}(s)}.
\end{equation*}

The \textit{weighted Laplacian} is defined by
\begin{equation*}
\Delta_f=\Delta+g(\nabla f,\nabla \cdot),
\end{equation*}
where $\Delta$ is the Laplacian defined as the minus of the divergence of the gradient.
The Laplacian comparison theorem in \cite{LMO:CompaFinsler} can be stated as follows (see \cite[Remark 3.10]{LMO:CompaFinsler}, and also \cite[Theorem 2.3]{KSa}):
\begin{thm}[\cite{LMO:CompaFinsler}]\label{thm:pointLapcomp}
We assume $\Ric_f^N\geq c^{-1}\,\kappa\, e^{-\frac{4(1-\varepsilon)f}{n-1}}g$.
Let $x\in \inte M$ and $v\in U_xM$.
Then for all $t \in ]0,\tau(v)[$ we have
\begin{equation}\label{eq:pointLapcomp}
\Delta_{f}\, \rho_{x}(\gamma_{v}(t))\geq H_{\kappa}(s_{f,v}(t))\,e^{-\frac{2(1-\eps)f(\gamma_{v}(t))}{n-1}}.
\end{equation}
\end{thm}

\subsection{Rigidity of Laplacian comparison from a single point}

The authors \cite{KSa} have shown the following (see \cite[Lemma 2.8]{KSa}):

\begin{lem}[\cite{KSa}]\label{lem:pointLaplacianRigidity}
Under the same setting as in Theorem \ref{thm:pointLapcomp},
assume that
the equality in \eqref{eq:pointLapcomp} holds at $t_{0}\in ]0, \tau(v)  [$.
Choose an orthonormal basis $\{e_{v,i}\}_{i=1}^{n}$ of $T_{x}M$ with $e_{v,n}=v$.
Let $\{Y_{v,i}\}^{n-1}_{i=1}$ and $\{E_{v,i}\}^{n-1}_{i=1}$ be the Jacobi fields and parallel vector fields along $\gamma_{v}$ with $Y_{v,i}(0)=0_x,\,Y_{v,i}'(0)=e_{v,i}$ and $E_{v,i}(0)=e_{v,i}$,
respectively.
Then the following properties hold on $[0,t_0]$:
\begin{enumerate}\setlength{\itemsep}{+0.7mm}
\item If $N=n$, then $f$ is constant, and
\begin{equation*}
Y_{v,i}(t)=\s_{\kappa\,e^{-\frac{4(1-\varepsilon)f}{n-1}}}(t)E_{v,i}(t);
\end{equation*} \label{enum:twisted curv}
\item if $N\neq 1,n$, then
\begin{equation*}
\eps=0,\quad f(\gamma_v(t))\equiv f(x),\quad Y_{v,i}(t)=\s_{\kappa \,e^{-\frac{4(1-\varepsilon)f(x)}{n-1}}}(t)E_{v,i}(t);
\end{equation*}\label{enum:relaxed twisted curv}
\item if $N=1$, then
\begin{equation*}
\eps=0,\quad Y_{v,i}(t)=\,\exp\left( \frac{f(\gamma_{v}(t))+f(x)}{n-1} \right)\,\s_{\kappa}(s_{f,v}(t))\,E_{v,i}(t).
\end{equation*} \label{enum:curv cond}
\end{enumerate}
\end{lem}

\subsection{Laplacian comparison from a single point with bounded density}

The authors \cite{KSa} have also proven the following (see \cite[Lemma 2.11]{KSa}, and also \cite[Theorem 3.9]{LMO:CompaFinsler}):
\begin{lem}[\cite{KSa}]\label{lem:finite pointed Laplacian comparison}
We assume
\begin{equation*}
\Ric_f^N\geq c^{-1}\,\kappa\, e^{-\frac{4(1-\varepsilon)f}{n-1}}g,\quad (1-\eps)f\leq (n-1)\delta
\end{equation*}
for $\delta \in \mathbb{R}$.
Let $x\in \inte M$ and $v\in U_xM$.
Then for all $t \in ]0,\tau(v)[$
we have
\begin{equation}\label{eq:finite pointed Laplacian comparison}
\Delta_{f}\, \rho_{x}(\gamma_{v}(t))\geq H_{\kappa}\left(e^{-2\delta}t\right)\,e^{-\frac{2(1-\eps)f(\gamma_{v}(t))}{n-1}}.
\end{equation}
\end{lem}

\begin{remark}\label{rem:finiteptrigid}
{\rm Assume that
the equality in $(\ref{eq:finite pointed Laplacian comparison})$ holds at $t_0$.
Then the equality in \eqref{eq:pointLapcomp} holds on $]0,t_{0}]$,
and $(1-\eps)f \circ \gamma_{v}=(n-1)\delta$ on $[0,t_{0}]$ (see \cite[Remark 2.12]{KSa}).}
\end{remark}

\section{Laplacian}\label{sec:Laplacian}

In this section,
we present a Laplacian comparison theorem for the distance function from the boundary,
which is a key ingredient of the proof of our main theorems.

\subsection{Riccati inequality}

Let $\rho_{\partial M}:M\to \mathbb{R}$ be the distance function from $\partial M$ defined as $\rho_{\partial M}(x):=d(x,\partial M)$.
For each $z\in \partial M$,
we denote by $\gamma_{z}:[0,T]\to M$ the unit speed geodesic with $\gamma_z(0)=z$ and $\dot{\gamma}_{z}(0)=u_z$.
We define a function $\tau:\partial M \to ]0,+\infty]$ as
\begin{equation*}
\tau(z):=\sup \{t>0 \mid \rho_{\partial M}(\gamma_{z}(t))=t \}.
\end{equation*}

We first show the following Riccati inequality:
\begin{lem}\label{lem:Riccati}
For all $t \in ]0,\tau(z)[$
we have
\begin{align}\label{eq:Riccati}
&\quad  \left(\bigl(e^{\frac{2(1-\eps)f}{n-1}}\,\Delta_{f}\rho_{\partial M}\bigl)(\gamma_{z}(t))\right)'\\ \notag
&\geq      e^{\frac{2(1-\eps)f(\gamma_{z}(t))}{n-1}}\,\Ric^{N}_{f}(\dot{\gamma}_{z}(t))+c\,e^{-\frac{2(1-\eps)f(\gamma_{z}(t))}{n-1}} \left(\bigl(e^{\frac{2(1-\eps)f}{n-1}}\,\Delta_{f}\rho_{\partial M}\bigl)(\gamma_{z}(t))   \right)^{2}.
\end{align}
\end{lem}
\begin{pf}
In the case of $N=n$,
the function $f$ is constant,
and the desired inequality is well-known.
In the case of $N=1$ with $\eps =0$,
\eqref{eq:Riccati} has been obtained in \cite{Sak} (see \cite[Lemma 3.1]{Sak}).
Hence we may assume $N\neq 1,n$.

Set $h_{f,v}:=\left(\Delta_{f}\rho_{\partial M}\right) \circ \gamma_{z}$ and $f_{z}:=f \circ \gamma_{z}$.
By applying the well-known Bochner formula to the distance function $\rho_{\partial M}$,
and by the Cauchy-Schwarz inequality,
\begin{align}\notag
0 &= \Ric^{\infty}_{f}(\dot{\gamma}_{z}(t))+\Vert \nabla^{2} \rho_{\partial M} \Vert^{2}\left(\gamma_{z}(t)\right)-g(\nabla \Delta_{f} \rho_{\partial M},\nabla \rho_{\partial M})(\gamma_z(t))\\  \label{eq:Cauchy-Schwarz}
&\geq \Ric^{N}_{f}(\dot{\gamma}_{z}(t))+\frac{f'_{z}(t)^{2}}{N-n}+\frac{\left(h_{f,z}(t)-f'_{z}(t)\right)^{2}}{n-1}-h'_{f,z}(t)\\  \notag
   &   =   \Ric^{N}_{f}(\dot{\gamma}_{z}(t))+ c\,h^{2}_{f,z}(t)-e^{-\frac{2(1-\eps)f_{z}(t)}{n-1}}\,\left(   e^{\frac{2(1-\eps)f_{z}(t)}{n-1}}  \,h_{f,z}(t)\right)' \\ \notag
   &\qquad +\frac{1}{n-1}\left( \sqrt{\frac{N-1}{N-n}}f'_{z}(t) -\eps \sqrt{\frac{N-n}{N-1}}h_{f,z}(t)  \right)^2\\ \label{eq:Bochner}
   &\geq \Ric^{N}_{f}(\dot{\gamma}_{z}(t))+ c\,h^{2}_{f,z}(t)-e^{-\frac{2(1-\eps)f_{z}(t)}{n-1}}\,\left(   e^{\frac{2(1-\eps)f_{z}(t)}{n-1}}  \,h_{f,z}(t)\right)'.
\end{align}
We arrive at the desired inequality \eqref{eq:Riccati}.
\end{pf}

\begin{remark}\label{rem:equalRiccati}
{\rm When $N\neq 1,n$,
we assume that
the equality in \eqref{eq:Riccati} holds at $t_{0}\in ]0,\tau(z)[$.
Then the equality in the Cauchy-Schwarz inequality in \eqref{eq:Cauchy-Schwarz} holds;
in particular,
\begin{equation*}
\nabla^2 \rho_{\partial M}=-\frac{\Delta \rho_{\partial M}}{n-1}\,g
\end{equation*}
on the orthogonal complement of $\nabla \rho_{\partial M}$ in $T_{\gamma_{z}(t_{0})}M$.
Moreover,
\begin{equation*}
\sqrt{\frac{N-1}{N-n}}f'_{z}(t_0) -\eps \sqrt{\frac{N-n}{N-1}}h_{f,z}(t_0)=0
\end{equation*}
since the equality in \eqref{eq:Bochner} holds.
In particular,
if $\eps =0$, then $f'_z(t_0)=0$,
and if $\eps \neq 0$,
\begin{equation*}
h_{f,z}(t_0)=\eps^{-1}\frac{N-1}{N-n}f'_z(t_0).
\end{equation*}}
\end{remark}

\subsection{Laplacian comparison theorem}
For $\kappa,\lambda \in \mathbb{R}$,
we denote by $\s_{\kappa,\lambda}(s)$ a unique solution to the Jacobi equation $\psi''(s)+\kappa\,\psi(s)=0$ with $\psi(0)=1$ and $\psi'(0)=-\lambda$.
We set
\begin{equation*}
C_{\kappa,\lambda}:=\inf\{s>0\mid \s_{\kappa,\lambda}(s)=0\}.
\end{equation*}
Notice that
$C_{\kappa,\lambda}$ is finite if and only if either
(1) $\kappa>0$; 
(2) $\kappa=0$ and $\lambda>0$;
or (3) $\kappa<0$ and $\lambda>\sqrt{\vert \kappa \vert}$,
and in this case,
we say that $\kappa$ and $\lambda$ satisfy the \textit{ball-condition}.
We also note that
they satisfy the ball-condition if and only if
there is a closed ball $B^{n}_{\kappa,\lambda}$ in the space form with constant curvature $\kappa$ whose boundary has constant mean curvature $(n-1)\lambda$.
The radius of $B^{n}_{\kappa,\lambda}$ is given by $C_{\kappa,\lambda}$.
We set
\begin{equation*}
H_{\kappa,\lambda}(s):=-c^{-1}\frac{\s'_{\kappa,\lambda}(s)}{\s_{\kappa,\lambda}(s)},
\end{equation*}
which enjoys the following Riccati equation:
\begin{equation}\label{eq:model Riccati}
H'_{\kappa,\lambda}(s)=c^{-1}\,\kappa+c\,H^2_{\kappa,\lambda}(s).
\end{equation}

Let us define functions $s_{f,z}:[0,\tau(z)]\to [0,\tau_{f}(z)]$ and $\tau_f:\partial M\to ]0,+\infty]$ by
\begin{equation*}
s_{f,z}(t):=\int^t_0\,e^{      -\frac{2(1-\eps)f(\gamma_z(\xi))}{n-1}       }\,\d\xi,  \quad \tau_{f}(z):=s_{f,z}(\tau(z)).
\end{equation*}
Let $t_{f,z}:[0,\tau_{f}(z)]\to [0,\tau(z)]$ be the inverse function of $t_{f,z}$.
Our Laplacian comparison theorem is stated as follows:
\begin{thm}\label{lem:Laplacian comparison}
Assume
\begin{equation}\label{eq:Lapcurvcond}
\Ric^{N}_{f}(\dot{\gamma}_{z}(t))\geq c^{-1}\,\kappa\,e^{   -\frac{4(1-\eps)f(\gamma_z(t))}{n-1}      },\quad H_{f,z}\geq c^{-1}\lambda e^{-\frac{2(1-\varepsilon)f(z)}{n-1}}
\end{equation}
for all $t \in ]0,\tau(z)[$.
Then for all $t \in ]0,\tau(z)[$ with $s_{f,z}(t) \in ]0,\min\{\tau_{f}(z),C_{\kappa,\lambda} \}[$,
\begin{equation}\label{eq:Laplacian comparison}
\Delta_{f}\rho_{\partial M}(\gamma_{z}(t)) \geq H_{\kappa,\lambda}(s_{f,z}(t))\,e^{      -\frac{2(1-\eps)f(\gamma_z(t))}{n-1}       }.
\end{equation}
\end{thm}
\begin{pf}
We define $F_{z}:]0,\tau(z)[ \to \mathbb{R}$ and $\hat{F}_{z}:]0,\tau_{f}(z)[ \to \mathbb{R}$ by 
\begin{equation*}
F_{z}:=\bigl(e^{\frac{2(1-\eps)f}{n-1}}\,\Delta_{f}\rho_{\partial M}\bigl) \circ \gamma_{z},\quad \hat{F}_{z}:=F_{z}\circ t_{f,z}.
\end{equation*}
From \eqref{eq:Riccati} and the curvature assumption,
for all $s \in ]0,\tau_{f}(z)[$,
\begin{align*}\label{eq:use of Riccati}
\hat{F}'_{z}(s)&    =   F'_{z}(t_{f,z}(s))\, e^{\frac{2(1-\eps)f\left(\gamma_{z}\left(    t_{f,z}(s)  \right)\right) }{n-1}}
\\ \notag
                      & \geq \Ric^{N}_{f}(\dot{\gamma}_{z}(t_{f,z}(s)))\,e^{\frac{4(1-\eps)f\left(\gamma_{z}\left(    t_{f,z}(s)  \right)\right) }{n-1}}
                      +c\,F^{2}_{z}(t_{f,z}(s))\\\notag
                      & \geq c^{-1}\kappa+c\,\hat{F}^{2}_{z}(s).
\end{align*}
The Riccati equation \eqref{eq:model Riccati} implies that
for all $s \in ]0,\min\{\tau_{f}(z) ,C_{\kappa,\lambda} \}  [$,
\begin{equation*}\label{eq:sharp Riccati}
\hat{F}'_{z}(s)-H'_{\kappa,\lambda}(s)\geq c\left(\hat{F}^{2}_{z}(s)-H^2_{\kappa,\lambda}(s) \right).
\end{equation*}

Let us consider a function $G_{\kappa,\lambda,z}:]0,\min\{\tau_{f}(z) ,C_{\kappa,\lambda} \}  [\to \mathbb{R}$ by
\begin{equation*}
G_{\kappa,\lambda,z}:=\mathfrak{s}^{2}_{\kappa,\lambda}\bigl( \hat{F}_{z}-H_{\kappa,\lambda} \bigl).
\end{equation*}
From (\ref{eq:sharp Riccati})
it follows that
\begin{align*}\label{eq:monotonicity}
G'_{\kappa,\lambda,z} &   =  2\,\mathfrak{s}_{\kappa,\lambda}\, \mathfrak{s}'_{\kappa,\lambda}\bigl( \hat{F}_{z}-H_{\kappa,\lambda} \bigl)+\mathfrak{s}^{2}_{\kappa,\lambda}\bigl( \hat{F}'_{z}-H'_{\kappa,\lambda} \bigl)\\ \notag
                                    &\geq 2\,\mathfrak{s}_{\kappa,\lambda}\, \mathfrak{s}'_{\kappa,\lambda}\bigl(\hat{F}_{z}-H_{\kappa,\lambda} \bigl)+c\,\mathfrak{s}^{2}_{\kappa,\lambda}\left(\hat{F}^{2}_{z}-H^{2}_{\kappa,\lambda}\right)\\ \notag
                                    &  =   c\,\mathfrak{s}^{2}_{\kappa,\lambda}\bigl(\hat{F}_{z}-H_{\kappa,\lambda} \bigl)^{2}\geq 0.   
\end{align*}
Since $G_{\kappa,\lambda,z}(s)\to e^{\frac{2(1-\eps)f(z)}{n-1}}\,H_{f,z}-c^{-1}\lambda$ as $s\to 0$,
the function $G_{\kappa,\lambda,z}$ is non-negative;
in particular,
$\hat{F}_{z} \geq H_{\kappa,\lambda}$ holds on $]0,\min\{\tau_{f}(z) ,C_{\kappa,\lambda} \}  [$.
This proves (\ref{eq:Laplacian comparison}).
\end{pf}

\begin{remark}\label{rem:equalLaplacian}
{\rm We assume that
the equality in $(\ref{eq:Laplacian comparison})$ holds at $t_{0}$.
Then $G_{\kappa,\lambda,z}(s_{f,z}(t_0))=0$.
From $G'_{\kappa,\lambda,z} \geq 0$
it follows that $G_{\kappa,\lambda,z}=0$ on $]0,s_{f,z}(t_0)]$;
in particular,
the equality in \eqref{eq:Riccati} holds on $]0,t_0]$ under \eqref{eq:Lapcurvcond} (see Remark \ref{rem:equalRiccati}).}
\end{remark}

\begin{remark}\label{rem:refLaplacian}
{\rm Theorem \ref{lem:Laplacian comparison} has been obtained by the second named author \cite{Sa2} and \cite{Sak} under \eqref{eq:Sa2setting} and \eqref{eq:Saksetting},
respectively (see \cite[Lemma 3.3]{Sa2}, \cite[Lemma 3.3]{Sak}).}
\end{remark}

Theorem \ref{lem:Laplacian comparison} leads us to the following:
\begin{lem}\label{lem:Cut point comparisons}
Let $\kappa$ and $\lambda$ satisfy the ball-condition.
Assume
\begin{equation*}
\Ric^{N}_{f}(\dot{\gamma}_{z}(t))\geq c^{-1}\,\kappa\,e^{   -\frac{4(1-\eps)f(\gamma_z(t))}{n-1}      },\quad H_{f,z}\geq c^{-1}\lambda e^{-\frac{2(1-\varepsilon)f(z)}{n-1}}
\end{equation*}
for all $t \in ]0,\tau(z)[$.
Then we have
\begin{equation*}
\tau_{f}(z) \leq C_{\kappa,\lambda}.
\end{equation*}
Moreover,
if $(1-\eps)f\circ \gamma_{z}\leq (n-1)\delta$ on $]0,\tau(z)[$ for $\delta \in \mathbb{R}$,
then
\begin{equation*}
\tau(z) \leq C_{\kappa\,e^{-4\delta},\lambda\,e^{-2\delta}}.
\end{equation*}
\end{lem}
\begin{pf}
One can prove it by using Theorem \ref{lem:Laplacian comparison} instead of \cite[Lemma 3.3]{Sak} along the line of the proof of \cite[Lemma 3.5]{Sak}.
We omit the proof. 
\end{pf}

\begin{remark}
{\rm Due to Lemma \ref{lem:Cut point comparisons},
one can drop the restriction $s_{f,z}(t) \in ]0,\min\{\tau_{f}(z),C_{\kappa,\lambda} \}[$ in Theorem \ref{lem:Laplacian comparison}.}
\end{remark}

\subsection{Rigidity of Laplacian comparison}

Let $A_{u_z}$ stand for the shape operator on $\partial M$ with respect to $u_z$.
We examine the equality case of Theorem \ref{lem:Laplacian comparison}.
\begin{lem}\label{lem:LaplacianRigidity}
Under the same setting as in Theorem \ref{lem:Laplacian comparison},
assume that
the equality in $(\ref{eq:Laplacian comparison})$ holds at $t_{0}\in ]0, \tau(z)  [$.
Choose an orthonormal basis $\{e_{z,i}\}_{i=1}^{n-1}$ of $T_{z}\partial M$.
Let $\{Y_{z,i}\}^{n-1}_{i=1}$ and $\{E_{z,i}\}^{n-1}_{i=1}$ be the Jacobi fields and parallel vector fields along $\gamma_{z}$ with $Y_{z,i}(0)=e_{z,i},\,Y_{z,i}'(0)=-A_{u_{z}}e_{z,i}$ and $E_{z,i}(0)=e_{z,i}$,
respectively.
Then the following hold on $[0,t_0]$:
\begin{enumerate}\setlength{\itemsep}{+0.7mm}
\item If $N=n$, then $f$ is constant, and
\begin{equation*}
Y_{z,i}(t)=\s_{\kappa\,e^{-\frac{4(1-\eps)f}{n-1}},\lambda\,e^{-\frac{2(1-\eps)f}{n-1}}}(t)E_{z,i}(t);
\end{equation*} \label{enum:twisted curv}
\item if $N\neq 1,n$, then
\begin{equation*}
f(\gamma_z(t))=f(z)-\eps \frac{N-n}{N-1}c^{-1}\log \mathfrak{s}_{\kappa,\lambda}(s_{f,z}(t)),\quad Y_{z,i}(t)=\s^{\frac{c^{-1}}{n-1} \left( 1-\eps \frac{N-n}{N-1}   \right) }_{\kappa,\lambda}(s_{f,z}(t)) E_{z,i}(t);
\end{equation*}\label{enum:relaxed twisted curv}
\item if $N=1$, then
\begin{equation*}
\eps=0,\quad Y_{z,i}(t)=\exp\left( \frac{f(\gamma_{z}(t))-f(z)}{n-1} \right)\,\s_{\kappa,\lambda}(s_{f,z}(t))\,E_{z,i}(t).
\end{equation*} \label{enum:curv cond}
\end{enumerate}
\end{lem}
\begin{pf}
If $N=n$,
then $f$ is constant by the definition, and the rigidity of Jacobi fields is well-known for $\eps=1$.
If $N=1$ with $\eps=0$,
then the desired assertion has been proved by the second named author \cite{Sak} (see \cite[Lemmas 3.8 and 3.9]{Sak}).
We may assume $N\neq 1,n$.

We first show the rigidity of $f$.
We set $h_{f,z}:=\left(\Delta_{f}\rho_{\partial M}\right) \circ \gamma_{z}$ and $f_{z}:=f \circ \gamma_{z}$.
Since the equality in $(\ref{eq:Laplacian comparison})$ holds at $t_{0}$,
the equalities in (\ref{eq:Riccati}) and (\ref{eq:Laplacian comparison}) also hold on $]0,t_{0}]$ (see Remarks \ref{rem:equalRiccati} and \ref{rem:equalLaplacian}).
If $\eps=0$,
then $f'_{z}(t)=0$ for each $t\in ]0,t_{0}]$,
and hence we already possess the desired formula.
Let us consider the case of $\eps \neq 0$.
Then $h_{f,z}(t)$ is equal to
\begin{equation}\label{eq:RigidweightLap}
\eps^{-1}\frac{N-1}{N-n}f'_z(t)=H_{\kappa,\lambda}(s_{f,z}(t))\,e^{      -\frac{2(1-\eps)f_z(t)}{n-1}       }.
\end{equation}
This implies
\begin{equation}\label{eq:densityrigid}
f_{z}(t)=f(z)-\eps \frac{N-n}{N-1}c^{-1}\log \mathfrak{s}_{\kappa,\lambda}(s_{f,z}(t)),
\end{equation}
which is the desired property.

We next investigate the rigidity of Jacobi fields.
We see
\begin{equation*}
\nabla^2 \rho_{\partial M}=-\frac{\Delta \rho_{\partial M}}{n-1}\,g
\end{equation*}
on the orthogonal complement of $\nabla \rho_{\partial M}$ in $T_{\gamma_{z}(t)}M$.
Set
\begin{equation*}
\varphi:=-\frac{\Delta \rho_{\partial M}}{n-1},\quad \varphi_{z}:=\varphi \circ \gamma_{z}.
\end{equation*}
The radial curvature equation (see e.g., \cite[Corollary 3.2.10]{Pet:RiemannianGeo}) yields
\begin{equation}\label{eq:radcurveq}
R(E_{z,i},\dot{\gamma}_z) \dot{\gamma}_z=-(\varphi'_{z}+\varphi^{2}_{z})E_{z,i}.
\end{equation}
On the other hand,
\begin{equation*}
\varphi_{z}(t)=\frac{1}{n-1}  \left(f'_{z}(t)-H_{\kappa,\lambda}(s_{f,z}(t)) \,e^{-\frac{2(1-\eps)f_z(t)}{n-1}}  \right)=\frac{1}{n-1}  \left(f_{z}(t)+\log \mathfrak{s}^{c^{-1}}_{\kappa,\lambda}(s_{f,z}(t)) \right)'
\end{equation*}
since the equality in \eqref{eq:Laplacian comparison} holds.
From direct computations and \eqref{eq:RigidweightLap},
we deduce
\begin{align}\label{hessian rigidity}
\varphi'_{z}(t)+\varphi^{2}_{z}(t)&=\frac{1}{n-1}\left( f''_z(t)-\frac{f'_{z}(t)^2}{N-n}-c^{-1}\,\kappa\,e^{-\frac{4(1-\eps)f_z(t)}{n-1}}    \right)\\ \notag
&\qquad +\frac{1}{n-1}\left( \sqrt{\frac{N-1}{N-n}}f'_z(t)-\eps \sqrt{\frac{N-n}{N-1}}H_{\kappa,\lambda}(s_{f,z}(t)) \,e^{-\frac{2(1-\eps)f_z(t)}{n-1}} \right)^2 \\ \notag
&=\frac{1}{n-1}\left( f''_z(t)-\frac{f'_{z}(t)^2}{N-n}-c^{-1}\,\kappa\,e^{-\frac{4(1-\eps)f_z(t)}{n-1}}    \right)=\frac{F''_{\kappa,\lambda,z}(t)}{F_{\kappa,\lambda,z}(t)},
\end{align}
where
\begin{equation}\label{eq:modelJacob}
F_{\kappa,\lambda,z}(t):=\exp \left(\frac{f_z(t)-f(z)}{n-1}\right)\,\s_{\kappa,\lambda}^{\frac{c^{-1}}{n-1}}(s_{f,z}(t)).
\end{equation}
The equality \eqref{hessian rigidity} together with \eqref{eq:radcurveq} tells us that $Y_{z,i}=F_{\kappa,\lambda,z}\,E_{z,i}$.
Substituting \eqref{eq:densityrigid} into \eqref{eq:modelJacob},
we arrive at
\begin{equation*}
F_{\kappa,\lambda,z}(t)=\s^{ \frac{c^{-1}}{n-1} \left( 1-\eps \frac{N-n}{N-1}   \right)  }_{\kappa,\lambda}(s_{f,z}(t)).
\end{equation*}
This completes the proof.
\end{pf}

\begin{remark}\label{rem:equality in finite Laplacian comparison}
{\rm Let us compare the rigidity phenomenon of Lemma \ref{lem:LaplacianRigidity} for $N\neq 1,n$ with that of Lemma \ref{lem:pointLaplacianRigidity}.
In such a case, the equality in \eqref{eq:pointLapcomp} holds only when $\eps=0$.
On the other hand,
it is possible that
the equality in (\ref{eq:Laplacian comparison}) holds for every $\eps$.}
\end{remark}

\subsection{Laplacian comparison with bounded density}\label{sec:Bounded cases}
We now investigate Laplacian comparisons under a boundedness of density.
Let us say that $\kappa$ and $\lambda$ satisfy the \textit{convex-ball-condition}
if they satisfy the ball-condition and $\lambda \geq 0$.
We say that
$\kappa$ and $\lambda$ satisfy the \textit{monotone-condition}
when $H_{\kappa,\lambda}\geq 0$ and $H'_{\kappa,\lambda}\geq 0$ on $[0,C_{\kappa,\lambda}[$.
We see that
they satisfy the monotone-condition
if and only if either
(1) they satisfy the convex-ball-condition;
or (2) $\kappa \leq 0$ and $\lambda=\sqrt{\vert \kappa \vert}$.
For $\kappa$ and $\lambda$ satisfying the monotone-condition,
if $\kappa = 0,\,\lambda=0$,
then $H_{\kappa,\lambda}=0$ on $[0,+\infty[$;
otherwise,
$H_{\kappa,\lambda} > 0$ on $]0,C_{\kappa,\lambda}[$.
We also say that
they satisfy the \textit{weakly-monotone-condition}
if $H'_{\kappa,\lambda} \geq 0$ on $[0,C_{\kappa,\lambda}[$.
Notice that
they satisfy the weakly-monotone-condition
if and only if either
(1) $\kappa \geq 0$;
or (2) $\kappa<0$ and $\vert \lambda \vert \geq \sqrt{\vert \kappa \vert}$.
In particular,
if $\kappa$ and $\lambda$ satisfy the ball-condition,
then they also satisfy the weakly-monotone-condition.
For $\kappa$ and $\lambda$ satisfying the weakly-monotone-condition,
if $\kappa \leq 0$ and $\vert \lambda \vert=\sqrt{\vert \kappa \vert}$,
then $H_{\kappa,\lambda}=(n-1)\lambda$ on $[0,+\infty[$;
otherwise,
$H'_{\kappa,\lambda} > 0$ on $[0,C_{\kappa,\lambda}[$.
From the same calculation as in the proof of \cite[Lemma 4.1]{Sak}, we can derive the following estimates from Theorem \ref{lem:Laplacian comparison}:
\begin{lem}\label{lem:finite Laplacian comparison}
Let $\kappa$ and $\lambda$ satisfy the weakly-monotone-condition.
We assume that
\begin{equation*}
\Ric^{N}_{f}(\gamma'_{z}(t))\geq (n-1)\kappa\, e^{-\frac{4 (1-\eps)f(\gamma_{z}(t))}{n-1}},\quad H_{f,z}\geq (n-1)\lambda e^{-\frac{2 (1-\eps)f(z)}{n-1}}
\end{equation*}
for all $t \in ]0,\tau(z)[$.
Suppose additionally that
$(1-\eps)f\circ \gamma_{z}\leq (n-1)\delta$ on $]0,\tau(z)[$ for $\delta \in \mathbb{R}$.
Then for all $t \in ]0,\tau(z)[$
we have
\begin{equation}\label{eq:finite Laplacian comparison}
\Delta_{f}\rho_{\partial M}(\gamma_{z}(t))  \geq H_{\kappa,\lambda}(e^{-2\delta}t)\,  e^{-\frac{2(1-\eps) f(\gamma_{z}(t))}{n-1}}.
\end{equation}
Moreover,
if $\kappa$ and $\lambda$ satisfy the monotone-condition,
then
\begin{equation}\label{eq:strictly finite Laplacian comparison}
\Delta_{f}\rho_{\partial M}(\gamma_{z}(t))  \geq H_{\kappa,\lambda}(e^{-2\delta}t)\,e^{-2\delta}.
\end{equation}
\end{lem}

\begin{remark}\label{rem:equality in finite Laplacian comparison}
{\rm Assume that
the equality in $(\ref{eq:finite Laplacian comparison})$ holds at $t_{0}$.
Then the equality in (\ref{eq:Laplacian comparison}) also holds (see Lemma \ref{lem:LaplacianRigidity}).
Moreover,
if either (1) $\kappa>0$; or (2) $\kappa \leq 0$ and $\vert \lambda \vert >\sqrt{\vert \kappa \vert}$,
then $(1-\eps)f\circ \gamma_{z}=(n-1)\delta$ on $[0,t_{0}]$ (cf. \cite[Remark 4.2]{Sak}).}
\end{remark}

\begin{remark}\label{rem:equality in strictly finite Laplacian comparison}
{\rm Assume that
the equality in $(\ref{eq:strictly finite Laplacian comparison})$ holds at $t_0$ under the monotone condition.
Then the equality in (\ref{eq:finite Laplacian comparison}) holds (see Remark \ref{rem:equality in finite Laplacian comparison}).
Moreover,
if either (1) $\kappa$ and $\lambda$ satisfy the convex-ball-condition; or (2) $\kappa<0$ and $\lambda=\sqrt{\vert \kappa \vert}$,
then $(1-\eps)(f\circ \gamma_{z})(t_{0})=(n-1)\delta$ (cf. \cite[Remark 4.3]{Sak}).}
\end{remark}

For $p\in ]1,+\infty[$,
the \textit{weighted $p$-Laplacian} is defined by
\begin{equation*}
\Delta_{f,p}:=-e^{f}\,\Div \,\left(e^{-f} \Vert \nabla \cdot \Vert^{p-2}\, \nabla \cdot \right).
\end{equation*}
By the same calculation and argument as in the proof of \cite[Lemma 4.4, Proposition 4.6]{Sak},
we have the following assertion:
\begin{prop}\label{prop:global finite p-Laplacian comparison}
Let $p\in ]1,+\infty[$.
Let $\kappa$ and $\lambda$ satisfy the monotone-condition.
Assume
\begin{align*}
\Ric_f^N\geq c^{-1}\,\kappa\, e^{-\frac{4(1-\varepsilon)f}{n-1}}g,\quad H_{f,\partial M}\geq c^{-1}\lambda e^{-\frac{2(1-\varepsilon)f}{n-1}}, \quad (1-\eps)f\leq (n-1)\delta
\end{align*}
for $\delta \in \mathbb{R}$.
Set $\rho_{\partial M,\delta}:=e^{-2\delta}\,\rho_{\partial M}$.
Let $\varphi:[0,+\infty[\to \mathbb{R}$ be a monotone increasing smooth function.
Then we have
\begin{equation*}
\Delta_{f,p}\, \left(\varphi \circ \rho_{\partial M,\delta} \right) \geq -e^{-2p\delta} \left[ \bigl( (\varphi' )^{p-1} \bigl)' -H_{\kappa,\lambda}\,  \left(\varphi' \right)^{p-1}\right]\circ \rho_{\partial M,\delta}
\end{equation*}
in the distribution sense.
\end{prop}

\subsection{Laplacian comparison with radial density}\label{sec:Radial cases}
Next,
we are concerned with the case where $f$ is $\partial M$-radial (i.e., there is a smooth function $\phi_{f}:[0,+\infty[\to \mathbb{R}$ such that $f=\phi_{f} \circ \rho_{\partial M}$).
Consider a function $\rho_{\partial M,f}:M\to \mathbb{R}$ defined by
\begin{equation*}\label{eq:weighted distance function from the boundary}
\rho_{\partial M,f}(x):=\inf_{z\in \partial M} \, \int^{\rho_{\partial M}(x)}_{0}\, e^{-\frac{2(1-\eps)f(\gamma_{z}(\xi))}{n-1}}\,\d \xi,
\end{equation*}
where the infimum is taken over all foot points $z\in \partial M$ of $x$ (i.e., $d(x,z)=\rho_{\partial M}(x)$).
A similar calculation to the proof of \cite[Lemma 4.9, Proposition 4.10]{Sak} yields the following inequality.
We omit the proof.
\begin{prop}\label{prop:global radial p-Laplacian comparison}
Let $p\in ]1,+\infty[$.
We assume
\begin{align*}
\Ric_f^N\geq c^{-1}\,\kappa\, e^{-\frac{4(1-\varepsilon)f}{n-1}}g,\quad H_{f,\partial M}\geq c^{-1}\lambda e^{-\frac{2(1-\varepsilon)f}{n-1}}.
\end{align*}
Suppose that
$f$ is $\partial M$-radial.
Let $\varphi:[0,+\infty[\to \mathbb{R}$ denote a monotone increasing smooth function.
Then we have
\begin{equation*}
\Delta_{\left( 1-\frac{2(p-1)(1-\eps)}{n-1}  \right)f,p}\,  \bigl(  \varphi \circ \rho_{\partial M,f} \bigl)  \geq -e^{\frac{-2 p(1-\eps) f}{n-1}} \left[ \bigl( (\varphi' )^{p-1} \bigl)'-H_{\kappa,\lambda} \,  (\varphi' )^{p-1}\right]\circ \rho_{\partial M,f}
\end{equation*}
in the distribution sense.
\end{prop}

\section{Splitting}\label{sec:Splitting theorems}
Hereafter,
we introduce our main theorems.
Once the Laplacian comparison is established,
we can prove the main results by using them along the line of the proof of the corresponding results in \cite{Sak} under \eqref{eq:Saksetting}.
We will just present their statements,
and outline or omit the proof.

In this section,
we show several splitting theorems.
Our first main result is the following splitting theorem,
which has been originally proven by Kasue \cite{K2}, Croke-Kleiner \cite{CK} in the unweighted case (see \cite[Theorem C]{K2}, \cite[Theorem 2]{CK}).
We will denote by $g_{\partial M}$ the induced Riemannian metric on $\partial M$.
\begin{thm}\label{thm:splitting theorem}
Let $\kappa \leq 0$ and $\lambda:=\sqrt{\vert \kappa \vert}$.
We assume
\begin{align*}
\Ric_f^N\geq c^{-1}\,\kappa\, e^{-\frac{4(1-\varepsilon)f}{n-1}}g,\quad H_{f,\partial M}\geq c^{-1}\lambda e^{-\frac{2(1-\varepsilon)f}{n-1}}.
\end{align*}
Suppose that
$(1-\eps)f$ is bounded from above.
If we have $\tau(z_{0})=\infty$ for some $z_{0}\in \partial M$,
then $M$ is diffeomorphic to $[0,+\infty[\times \partial M$,
and the following properties hold:
\begin{enumerate}\setlength{\itemsep}{+0.7mm}
\item If $N=n$, then $f$ is constant, and
\begin{equation*}
g=dt^2+\s^2_{\kappa\,e^{-\frac{4(1-\eps)f}{n-1}},\lambda\,e^{-\frac{2(1-\eps)f}{n-1}}}(t)g_{\partial M};
\end{equation*} \label{enum:twisted curv}
\item if $N\neq 1,n$, then for any $z\in \partial M$
\begin{equation*}
f(\gamma_z(t))=f(z)-\eps \frac{N-n}{N-1}c^{-1}\log \mathfrak{s}_{\kappa,\lambda}(s_{f,z}(t)),\quad g=dt^2+\s^{2\frac{c^{-1}}{n-1} \left( 1-\eps \frac{N-n}{N-1}   \right) }_{\kappa,\lambda}(s_{f,z}(t))g_{\partial M};
\end{equation*}\label{enum:relaxed twisted curv}
\item if $N=1$, then
\begin{equation*}
\eps=0,\quad g=dt^2+\exp\left(2 \frac{f(\gamma_{z}(t))-f(z)}{n-1} \right)\,\s^2_{\kappa,\lambda}(s_{f,z}(t))\,g_{\partial M}.
\end{equation*} \label{enum:curv cond}
\end{enumerate}
\end{thm}
\begin{pf}
We only sketch the proof.
If $N=n$,
then the desired statement follows from \cite{K2}, \cite{CK}.
If $N=1$,
then it has been obtained by the second named author \cite{Sak} (see \cite[Theorem 1.1]{Sak}).
Thus it is enough to discuss the case of $N\neq 1,n$.

The proof is similar to that of $N=1$ in \cite{Sak}.
For the statement that $M$ is diffeomorphic to $[0,+\infty[\times \partial M$,
we can prove it only by replacing the role of \cite[Lemmas 2.9 and 3.3]{Sak} with Lemma \ref{lem:finite pointed Laplacian comparison} and Theorem \ref{lem:Laplacian comparison} along the line of the proof of \cite[Theorem 1.1]{Sak}.
Then the equality in Theorem \ref{lem:Laplacian comparison} holds over $M$.
Lemma \ref{lem:LaplacianRigidity} \ref{enum:relaxed twisted curv} leads us to the conclusion.
\end{pf}

For $\kappa>0$ and $\lambda<0$,
we set
\begin{equation*}
D_{\kappa,\lambda}:=\inf\, \{\,s>0 \mid  \mathfrak{s}'_{\kappa,\lambda}(s)=0\,\}.
\end{equation*}
We also obtain the following splitting theorem,
which has been formulated by Kasue \cite{K2} in the unweighted setting (see \cite[Theorem B]{K2}):
\begin{thm}\label{thm:disconnected splitting2}
Let $\kappa>0$.
We assume
\begin{align*}
\Ric_f^N\geq c^{-1}\,\kappa\, e^{-\frac{4(1-\varepsilon)f}{n-1}}g,\quad H_{f,\partial M}\geq c^{-1}\lambda e^{-\frac{2(1-\varepsilon)f}{n-1}},\quad (1-\eps)f\leq (n-1)\delta
\end{align*}
for $\delta \in \mathbb{R}$.
Let $\partial M$ be disconnected,
and let $\{\partial M_{i}\}_{i=1,2,\dots}$ denote the connected components of $\partial M$.
Let $\partial M_{1}$ be compact.
Then $\lambda<0$, and
\begin{equation*}
\inf_{i=2,3,\dots}\, d(\partial M_{1},\partial M_{i}) \leq 2D_{\kappa\,e^{-4\delta},\lambda\,e^{-2\delta}}.
\end{equation*}
Moreover,
if the equality holds,
then $M$ is diffeomorphic to $[0,2D_{\kappa\,e^{-4\delta},\lambda\,e^{-2\delta}}]\times \partial M_1$,
and
\begin{equation*}
(1-\eps)f=(n-1)\delta,\quad g=dt^2+\s^2_{\kappa\,e^{-4\delta},\lambda\,e^{-2\delta}}(t)g_{\partial M_1}.
\end{equation*}
\end{thm}
\begin{pf}
The second named author \cite{Sak} has proved this result when $N=1$ (see \cite[Theorem 5.8]{Sak}).
Along the line of its proof,
we can prove the desired assertion by using \eqref{eq:finite Laplacian comparison} instead of \cite[(4.1)]{Sak} (see also Remark \ref{rem:equality in finite Laplacian comparison}).
\end{pf}

\begin{remark}
{\rm The second named author \cite{Sa2} has shown Theorems \ref{thm:splitting theorem} and \ref{thm:disconnected splitting2} under the condition \eqref{eq:Sa2setting} with $N\in [n,+\infty[$ and $\eps=1$ (see \cite[Theorems 1.4 and 6.14]{Sa2}).}
\end{remark}

\section{Inscribed radius}\label{sec:Inscribed radii}
We next study the inscribed radius $\IR M$ which is defined as the supremum of the distance function $\rho_{\partial M}$ from the boundary over $M$.
Let us consider a conformally deformed Riemannian metric
\begin{equation*}
g_{f}:=e^{-\frac{4(1-\eps)f}{n-1}}g.
\end{equation*}
Let $\rho^{g_{f}}_{\partial M}$ and $\IR_{g_{f}} M$ stand for the distance function from the boundary and the inscribed radius induced from $g_{f}$,
respectively.
\begin{thm}\label{thm:inscribed radius rigidity}
Let $\kappa$ and $\lambda$ satisfy the ball-condition.
We assume
\begin{align*}
\Ric_f^N\geq c^{-1}\,\kappa\, e^{-\frac{4(1-\varepsilon)f}{n-1}}g,\quad H_{f,\partial M}\geq c^{-1}\lambda e^{-\frac{2(1-\varepsilon)f}{n-1}}.
\end{align*}
Then we have
\begin{equation*}\label{eq:inscribed radius rigidity}
\IR_{g_{f}} M \leq C_{\kappa,\lambda}.
\end{equation*}
If $\rho^{g_{f}}_{\partial M}(x_{0})=C_{\kappa,\lambda}$ for some $x_{0}\in M$,
then $M$ is diffeomorphic to a closed ball centered at $x_0$,
and the following hold:
\begin{enumerate}\setlength{\itemsep}{+0.7mm}
\item If $N=n$, then $f$ is constant, and
\begin{equation*}
g=dt^2+\s^2_{\kappa\,e^{-\frac{4(1-\varepsilon)f}{n-1}}}(t)g_{\mathbb{S}^{n-1}};
\end{equation*} \label{enum:twisted curv}
\item if $N\neq 1,n$, then $f$ is constant, and 
\begin{equation*}
\eps=0,\quad g=dt^2+\s^2_{\kappa \,e^{-\frac{4f}{n-1}}}(t)g_{\mathbb{S}^{n-1}};
\end{equation*}\label{enum:relaxed twisted curv}
\item if $N=1$, then $f$ is radial with respect to $x_0$, and
\begin{equation*}
\eps=0,\quad g=dt^2+\exp\left( 2\frac{f(\gamma_{v}(t))+f(x_0)}{n-1} \right)\,\s^2_{\kappa}(s_{f,v}(t))g_{\mathbb{S}^{n-1}},
\end{equation*}
here $\gamma_{v}:[0,\rho_{\partial M}(x_0)]\to M$ denotes the geodesic with $\gamma_{v}(0)=x_0$ and ${\dot\gamma}_{v}(0)=v$. 
\end{enumerate}
\end{thm}
\begin{pf}
Let us sketch the proof.
This is a weighted version of the result by Kasue \cite{K2} (see \cite[Theorem A]{K2}).
The result for $N=n$ can be directly derived from \cite{K2}.
Furthermore,
the second named author \cite{Sak} has shown Theorem \ref{thm:inscribed radius rigidity} when $N=1$ (see \cite[Theorem 1.2]{Sak}).
We may assume $N\neq 1,n$.

Similarly to Theorem \ref{thm:splitting theorem},
we can refer to the argument of \cite{Sak}.
In the proof of the claim that $M$ is diffeomorphic to a closed ball,
the key point is to show the subharmonicity of the function $\rho_{x_0}+\rho_{\partial M}$.
This is done by replacing the role of \cite[Lemmas 2.8 and 3.3]{Sak} with Theorems \ref{thm:pointLapcomp} and \ref{lem:Laplacian comparison} along the line of the proof of \cite[Theorem 1.2]{Sak}.
Then the equality in Theorem \ref{thm:pointLapcomp} holds on $M$,
and Lemma \ref{lem:pointLaplacianRigidity} implies the desired conclusion.
\end{pf}

\begin{remark}
{\rm Under the condition \eqref{eq:Sa2setting} with $N\in [n,+\infty[$ and $\eps=1$,
Li-Wei \cite{LW2} have proved Theorem \ref{thm:inscribed radius rigidity} when $\kappa=0$ (see \cite[Theorem 4]{LW2}).
In \cite{LW1},
they also have done when $\kappa<0$ (see \cite[Theorem 1.2]{LW1}). 
Finally,
the second named author \cite{Sa2} has done for all $\kappa$ and $\lambda$ satisfying the ball-condition (see \cite[Theorems 1.1]{Sa2}).
We also refer to \cite{BKMW} on the work in non-smooth setting.}
\end{remark}

Under a boundedness of $f$,
we conclude the following:
\begin{thm}\label{thm:finite inscribed radius rigidity}
Let $\kappa$ and $\lambda$ satisfy the ball-condition.
We assume
\begin{align*}
\Ric_f^N\geq c^{-1}\,\kappa\, e^{-\frac{4(1-\varepsilon)f}{n-1}}g,\quad H_{f,\partial M}\geq c^{-1}\lambda e^{-\frac{2(1-\varepsilon)f}{n-1}},\quad (1-\eps)f\leq (n-1)\delta
\end{align*}
for $\delta \in \mathbb{R}$.
Then we have
\begin{equation*}
\IR M \leq C_{\kappa\,e^{-4\delta},\lambda\,e^{-2\delta}}.
\end{equation*}
If we have $\rho_{\partial M}(x_{0})=C_{\kappa\,e^{-4\delta},\lambda\,e^{-2\delta}}$ for some $x_{0}\in M$,
then $(1-\eps)f=(n-1)\delta$ on $M$,
and the following properties hold:
\begin{enumerate}\setlength{\itemsep}{+0.7mm}
\item If $N=n$, then $M$ is isometric to $B^{n}_{\kappa\,e^{-4\delta},\lambda\,e^{-2\delta}}$;
\item if $N\neq n$, then $\eps =0$, and $M$ is isometric to $B^{n}_{\kappa\,e^{-4\delta},\lambda\,e^{-2\delta}}$.
\end{enumerate}
\end{thm}
\begin{pf}
The second named author \cite{Sak} has proved this result when $N=1$ (see \cite[Theorem 6.4]{Sak}).
Along the line of its proof,
one can prove the claim by using Lemma \ref{lem:finite pointed Laplacian comparison} and \eqref{eq:finite Laplacian comparison} instead of \cite[Lemma 2.9 and (4.1)]{Sak} (see also Remark \ref{rem:finiteptrigid}).
\end{pf}

\section{Volume}\label{sec:Volume growths}
This section is devoted to the study of volume comparisons.

\subsection{Volume elements}\label{sec:volume elements}
For $z\in \partial M$ and $t \in ]0,\tau(z)[$,
let $\theta(t,z)$ be the volume element of the $t$-level surface of $\rho_{\partial M}$ at $\gamma_{z}(t)$.
For $s \in ]0,\tau_{f}(z)[$,
we set
\begin{equation*}\label{eq:volume element}
\theta_{f}(t,z):=e^{-f(\gamma_{z}(t))}\, \theta(t,z),\quad \hat{\theta}_{f}(s,z):=\theta_{f}(t_{f,z}(s),z).
\end{equation*}
By Theorem \ref{lem:Laplacian comparison} and the same calculation as in \cite[Lemma 7.1]{Sak},
we obtain:
\begin{lem}\label{lem:volume element comparison}
Assume
\begin{equation*}
\Ric^{N}_{f}(\dot{\gamma}_{z}(t))\geq c^{-1}\,\kappa\,e^{   -\frac{4(1-\eps)f(\gamma_z(t))}{n-1}      },\quad H_{f,z}\geq c^{-1}\lambda e^{-\frac{2(1-\varepsilon)f(z)}{n-1}}
\end{equation*}
for all $t \in ]0,\tau(z)[$.
Then for all $s_{1},s_{2} \in [0,\tau_{f}(z)[$ with $s_{1}\leq s_{2}$
\begin{equation*}
\frac{\hat{\theta}_{f}(s_{2},z)}{ \hat{\theta}_{f}(s_{1},z)}\leq \frac{\mathfrak{s}^{c^{-1}}_{\kappa,\lambda}(s_{2})}{\mathfrak{s}^{c^{-1}}_{\kappa,\lambda}(s_{1})}.
\end{equation*}
In particular,
for all $s\in [0,\tau_{f}(z)[$ we have
\begin{equation*}
\hat{\theta}_{f}(s,z)\leq e^{-f(z)}\,\mathfrak{s}^{c^{-1}}_{\kappa,\lambda}(s).
\end{equation*}
\end{lem}

Furthermore,
in virtue of \eqref{eq:strictly finite Laplacian comparison},
we have the following (cf. \cite[Lemma 7.3]{Sak}):
\begin{lem}\label{lem:finite volume element comparison}
Let $\kappa$ and $\lambda$ satisfy the monotone-condition.
Assume
\begin{equation*}
\Ric^{N}_{f}(\dot{\gamma}_{z}(t))\geq c^{-1}\,\kappa\,e^{   -\frac{4(1-\eps)f(\gamma_z(t))}{n-1}      },\quad H_{f,z}\geq c^{-1}\lambda e^{-\frac{2(1-\varepsilon)f(z)}{n-1}},\quad (1-\eps)f \circ \gamma_{z} \leq (n-1)\delta
\end{equation*}
on $(0,\tau(z))$.
Then for all $t_{1},t_{2} \in [0,\tau(z)[$ with $t_{1}\leq t_{2}$
we have
\begin{equation*}\label{eq:Jacobi comparison}
\frac{\theta_{f}(t_{2},z)}{ \theta_{f}(t_{1},z)}\leq \frac{\mathfrak{s}^{c^{-1}}_{\kappa\,e^{-4\delta},\lambda\,e^{-2\delta}}(t_{2})}{\mathfrak{s}^{c^{-1}}_{\kappa\,e^{-4\delta},\lambda\,e^{-2\delta}}(t_{1})}.
\end{equation*}
In particular,
for all $t\in [0,\tau(z)[$
we have
\begin{equation*}
\theta_{f}(t,z)\leq e^{-f(z)}\,\mathfrak{s}^{c^{-1}}_{\kappa\,e^{-4\delta},\lambda\,e^{-2\delta}}(t).
\end{equation*}
\end{lem}

\subsection{Volume comparisons}\label{sec:Absolute volume comparisons}
We define
\begin{equation*}
m_{f}:=e^{-f}v_g,\quad m_{f,\partial M}:=e^{-f|_{\partial M}}v_{g_{\partial M}},
\end{equation*}
where $v_{g}$ and $v_{g_{\partial M}}$ are the Riemannian volume measure on $M$ and $\partial M$ induced from $g$ and $g_{\partial M}$,
respectively.
For $r>0$ we set
\begin{equation*}
B_{r}(\partial M):=\{\,x\in M \mid \rho_{\partial M}(x) \leq r \,\},\quad B^{f}_{r}(\partial M):=\{\,x\in M \mid \rho_{\partial M,f}(x) \leq r \,\}.
\end{equation*}
We notice the following (cf. \cite[Lemma 5.1]{Sa2}, \cite[Lemma 7.5]{Sak}):
If $\partial M$ is compact, then
\begin{align*}
m_{\left( 1+\frac{2(1-\eps)}{n-1}  \right)f}\left(B^{f}_{r}(\partial M)\right)&=\int_{\partial M}\,\int^{\min\{r,\tau_f(z)\}}_{0}\,\hat{\theta}_{f}(s,z)\,\d s\, \d v_{g_{\partial M}},\\
m_{f}\left(B_{r}(\partial M)\right)&=\int_{\partial M}\,\int^{\min\{r,\tau(z)\}}_{0}\,\theta_{f}(t,z)\,\d t\, \d v_{g_{\partial M}}.
\end{align*}
We also set
\begin{equation*}
\mathcal{S}_{\kappa,\lambda}(r):=\int^{\min\{r,C_{\kappa,\lambda}\}}_{0}\,\mathfrak{s}^{c^{-1}}_{\kappa,\lambda}(\xi)\,\d \xi.
\end{equation*}

Lemmas \ref{lem:volume element comparison}, \ref{lem:finite volume element comparison} together with the same argument as in \cite[Subsection 7.2]{Sak} for $N=1$ tell us the following Heintze-Karcher type comparisons (cf. \cite[Lemmas 7.6 and 7.7]{Sak}):
\begin{prop}\label{lem:absolute volume comparison}
Assume
\begin{align*}
\Ric_f^N\geq c^{-1}\,\kappa\, e^{-\frac{4(1-\varepsilon)f}{n-1}}g,\quad H_{f,\partial M}\geq c^{-1}\lambda e^{-\frac{2(1-\varepsilon)f}{n-1}}.
\end{align*}
Let $\partial M$ be compact.
Then for all $r>0$
we have
\begin{equation*}\label{eq:absolute volume comparison}
m_{\left( 1+\frac{2(1-\eps)}{n-1}  \right)f}\left(B^{f}_{r}(\partial M)\right) \leq \mathcal{S}_{\kappa,\lambda}(r)\,m_{f,\partial M}(\partial M).
\end{equation*}
\end{prop}

\begin{prop}\label{lem:finite absolute volume comparison}
Let $\kappa$ and $\lambda$ satisfy the monotone-condition.
Assume
\begin{align*}
\Ric_f^N\geq c^{-1}\,\kappa\, e^{-\frac{4(1-\varepsilon)f}{n-1}}g,\quad H_{f,\partial M}\geq c^{-1}\lambda e^{-\frac{2(1-\varepsilon)f}{n-1}},\quad (1-\eps)f\leq (n-1)\delta
\end{align*}
for $\delta \in \mathbb{R}$.
Let $\partial M$ be compact.
Then for all $r>0$
we have
\begin{equation*}\label{eq:finite absolute volume comparison}
m_{f}(B_{r}(\partial M)) \leq \mathcal{S}_{\kappa\,e^{-4\,\delta},\lambda\,e^{-2\,\delta}}(r)\,m_{f,\partial M}(\partial M).
\end{equation*}
\end{prop}

\begin{remark}
{\rm Under the curvature condition \eqref{eq:Sa2setting},
Bayle \cite{B} has stated Proposition \ref{lem:absolute volume comparison} without proof (see \cite[Theorem E.2.2]{B}, and also \cite{M}, \cite{Mo}).}
\end{remark}

Also,
Lemmas \ref{lem:volume element comparison} and \ref{lem:finite volume element comparison} together with the same argument as in \cite[Subsection 7.3]{Sak} yield the following relative volume comparisons (cf. \cite[Theorems 1.3 and 7.9]{Sak}):
\begin{prop}\label{thm:relative volume comparison}
Under the same setting as in Proposition \ref{lem:absolute volume comparison},
for all $r,R>0$ with $r\leq R$
we have
\begin{equation*}\label{eq:relative volume comparison}
\frac{m_{\left( 1+\frac{2(1-\eps)}{n-1}  \right)f}(B^{f}_{R}(\partial M))}{m_{\left( 1+\frac{2(1-\eps)}{n-1}  \right)f}(B^{f}_{r}(\partial M))} \leq \frac{\mathcal{S}_{\kappa,\lambda}(R)}{\mathcal{S}_{\kappa,\lambda}(r)}.
\end{equation*}
\end{prop}

\begin{prop}\label{thm:finite relative volume comparison}
Under the same setting as in Proposition \ref{lem:finite absolute volume comparison},
for all $r,R>0$ with $r\leq R$ we have
\begin{equation*}\label{eq:finite relative volume comparison}
\frac{m_{f}(B_{R}(\partial M))}{m_{f}(B_{r}(\partial M))} \leq \frac{\mathcal{S}_{\kappa\,e^{-4\delta},\lambda\,e^{-2\delta}}(R)}{\mathcal{S}_{\kappa\,e^{-4\delta},\lambda\,e^{-2\delta}}(r)}.
\end{equation*}
\end{prop}

\subsection{Volume rigidity}\label{sec:Volume growth rigidity}
By the same method as in the proof of \cite[Theorem 7.11]{Sak} concerning $N=1$,
we see the following rigidity result for the equality case of Propositions \ref{lem:absolute volume comparison} and \ref{thm:relative volume comparison}:
\begin{thm}\label{thm:volume growth rigidity}
Under the same setting as in Propositions \ref{lem:absolute volume comparison} and \ref{thm:relative volume comparison},
if $\kappa$ and $\lambda$ do not satisfy the ball-condition,
and if
\begin{equation}\label{eq:assumption of volume growth rigidity}
\liminf_{r\to \infty}\frac{m_{\left( 1+\frac{2(1-\eps)}{n-1}  \right)f}(B^{f}_{r}(\partial M))}{\mathcal{S}_{\kappa,\lambda}(r)}\geq m_{f,\partial M}(\partial M),
\end{equation}
then $M$ is diffeomorphic to $[0,+\infty[\times \partial M$,
and the following hold:
\begin{enumerate}\setlength{\itemsep}{+0.7mm}
\item If $N=n$, then $f$ is constant, and
\begin{equation*}
g=dt^2+\s^2_{\kappa\,e^{-\frac{4(1-\eps)f}{n-1}},\lambda\,e^{-\frac{2(1-\eps)f}{n-1}}}(t)g_{\partial M};
\end{equation*} \label{enum:twisted curv}
\item if $N\neq 1,n$, then for any $z\in \partial M$
\begin{equation*}
f(\gamma_z(t))=f(z)-\eps \frac{N-n}{N-1}c^{-1}\log \mathfrak{s}_{\kappa,\lambda}(s_{f,z}(t)),\quad g=dt^2+\s^{2\frac{c^{-1}}{n-1} \left( 1-\eps \frac{N-n}{N-1}   \right) }_{\kappa,\lambda}(s_{f,z}(t))g_{\partial M};
\end{equation*}\label{enum:relaxed twisted curv}
\item if $N=1$, then
\begin{equation*}
\eps=0,\quad g=dt^2+\exp\left(2 \frac{f(\gamma_{z}(t))-f(z)}{n-1} \right)\,\s^2_{\kappa,\lambda}(s_{f,z}(t))\,g_{\partial M}.
\end{equation*} \label{enum:curv cond}
\end{enumerate}
\end{thm}

Thanks to Theorems \ref{thm:splitting theorem} and \ref{thm:finite inscribed radius rigidity},
the same argument as in the proof of \cite[Theorem 7.13]{Sak} leads to the following result for the equality case of Propositions \ref{lem:finite absolute volume comparison} and \ref{thm:finite relative volume comparison}:
\begin{thm}\label{thm:finite volume growth rigidity}
Under the same setting as in Propositions \ref{lem:finite absolute volume comparison} and \ref{thm:finite relative volume comparison},
if
\begin{equation*}\label{eq:assumption of finite volume growth rigidity}
\liminf_{r\to \infty}\frac{m_{f}(B_{r}(\partial M))}{\mathcal{S}_{\kappa\,e^{-4\,\delta},\lambda\,e^{-2\,\delta}}(r)}\geq m_{f,\partial M}(\partial M),
\end{equation*}
then the following hold:
\begin{enumerate}
\item if $\kappa$ and $\lambda$ satisfy the convex-ball-condition,
         then $M$ is isometric to $B^{n}_{\kappa e^{-4\delta},\lambda e^{-2\delta}}$,
         and $(1-\eps)f=(n-1)\delta$ on $M$;
\item if $\kappa \leq 0$ and $\lambda=\sqrt{\vert \kappa \vert}$,
         then $M$ is diffeomorphic to $[0,+\infty[\times \partial M$,
and the following rigidity properties hold:
\begin{enumerate}\setlength{\itemsep}{+0.7mm}
\item If $N=n$, then $f$ is constant, and
\begin{equation*}
g=dt^2+\s^2_{\kappa\,e^{-\frac{4(1-\eps)f}{n-1}},\lambda\,e^{-\frac{2(1-\eps)f}{n-1}}}(t)g_{\partial M};
\end{equation*} \label{enum:twisted curv}
\item if $N\neq 1,n$, then for any $z\in \partial M$
\begin{align*}
f(\gamma_z(t))&=f(z)-\eps \frac{N-n}{N-1}c^{-1}\log \mathfrak{s}_{\kappa,\lambda}(s_{f,z}(t)),\\
g&=dt^2+\s^{2\frac{c^{-1}}{n-1} \left( 1-\eps \frac{N-n}{N-1}   \right) }_{\kappa,\lambda}(s_{f,z}(t))g_{\partial M};
\end{align*}\label{enum:relaxed twisted curv}
\item if $N=1$, then
\begin{equation*}
\eps=0,\quad g=dt^2+\exp\left(2 \frac{f(\gamma_{z}(t))-f(z)}{n-1} \right)\,\s^2_{\kappa,\lambda}(s_{f,z}(t))\,g_{\partial M};
\end{equation*} \label{enum:curv cond}
\end{enumerate}
moreover, if $\kappa<0$, then $(1-\eps)f=(n-1)\delta$ on $M$.
\end{enumerate}
\end{thm}

\section{Spectrum}\label{sec:Spectrum}
In this last section,
we will collect comparison geometric results on the smallest Dirichlet eigenvalue for the weighted $p$-Laplacian $\Delta_{f,p}$.

\subsection{Eigenvalue rigidity}
Let $p\in ]1,+\infty[$.
A real number $\nu$ is said to be a \textit{Dirichlet eigenvalue of $\Delta_{f,p}$}
if there exists $\psi \in W^{1,p}_{0}(M,m_{f}) \setminus \{0\}$
such that $\Delta_{f,p} \psi=\nu \vert \psi \vert^{p-2}\,\psi$ in the distribution sense.
For $\psi \in W^{1,p}_{0}(M,m_{f})\setminus \{0\}$,
the associated \textit{Rayleigh quotient} is defined as
\begin{equation*}
R_{f,p}(\psi):=\frac{\int_{M}\, \Vert \nabla \psi \Vert^{p}\,\d\,m_{f}}{\int_{M}\,  \vert \psi \vert^{p}\,\d\,m_{f}}.
\end{equation*}
We set
\begin{equation*}
\nu_{f,p}(M):=\inf_{\psi} R_{f,p}(\psi),
\end{equation*}
where the infimum is taken over all $\psi \in W^{1,p}_{0}(M,m_{f}) \setminus \{0\}$.
If $M$ is compact,
then $\nu_{f,p}(M)$ coincides with the infimum of the set of all Dirichlet eigenvalues.

For $D\in ]0,C_{\kappa,\lambda}]\setminus \{+\infty\}$,
let $\nu_{p,\kappa,\lambda,D}$ be the positive minimum real number $\nu$ such that
there is a non-zero function $\varphi:[0,D]\to \mathbb{R}$ satisfying
\begin{align*}\label{eq:model space eigenvalue problem}
\left(\vert \varphi'(s)\vert^{p-2} \varphi'(s)\right)'&+(n-1)\frac{\mathfrak{s}'_{\kappa,\lambda}(s)}{\mathfrak{s}_{\kappa,\lambda}(s)}\,\left(\vert \varphi'(s) \vert^{p-2} \varphi'(s)\right)\\ \notag
                                                                         &+\nu\, \vert \varphi(s)\vert^{p-2}\varphi(s)=0, \,\, \varphi(0)=0,\,\, \varphi'(D)=0.                           
\end{align*}

We recall the notion of the model spaces introduced by Kasue \cite{K3}.
We say that
$\kappa$ and $\lambda$ satisfy the \textit{model-condition}
if the equation $\mathfrak{s}'_{\kappa,\lambda}(s)=0$ has a positive solution.
Note that
$\kappa$ and $\lambda$ satisfy the model-condition if and only if either
(1) $\kappa>0$ and $\lambda<0$;
(2) $\kappa=0$ and $\lambda=0$;
or (3) $\kappa<0$ and $\lambda \in ]0,\sqrt{\vert \kappa \vert}[$.
Let $\kappa$ and $\lambda$ satisfy the ball-condition or the model-condition,
and let $M$ be compact.
When they satisfy the model-condition,
we set $D_{\kappa,\lambda}(M)$ as follows:
If $\kappa=0$ and $\lambda=0$,
then $D_{\kappa,\lambda}(M):=\IR M$;
otherwise,
$D_{\kappa,\lambda}(M):=\inf \{s>0\mid  \mathfrak{s}'_{\kappa,\lambda}(s)=0\}$.
We say that $M$ is a \textit{$\left(\kappa,\lambda \right)$-equational model space}
if $M$ is isometric to either
(1) $B^n_{\kappa,\lambda}$ for $\kappa$ and $\lambda$ satisfying the ball-condition;
(2) $\left([0,2D_{\kappa,\lambda}(M)] \times \partial M_{1},ds^{2}+ \mathfrak{s}^{2}_{\kappa,\lambda}(s)\,g_{\partial M_1}\right)$ for $\kappa$ and $\lambda$ satisfying the model-condition,
and for some connected component $\partial M_{1}$ of $\partial M$;
or (3) the quotient space $\left([0,2D_{\kappa,\lambda}(M)] \times \partial M,ds^{2}+ \mathfrak{s}^{2}_{\kappa,\lambda}(s)\,g_{\partial M}\right)/G_{\sigma}$ for $\kappa$ and $\lambda$ satisfying the model-condition,
and for some involutive isometry $\sigma$ of $\partial M$ without fixed points,
where $G_{\sigma}$ denotes the isometry group on $\left([0,2D_{\kappa,\lambda}(M)] \times \partial M,ds^{2}+ \mathfrak{s}^{2}_{\kappa,\lambda}(s)\,g_{\partial M}\right)$
whose elements consist of identity and the involute isometry $\hat{\sigma}$ defined by $\hat{\sigma}(s,z):=\left(2D_{\kappa,\lambda}(M)-s,\sigma(z)\right)$.

We now set
\begin{equation*}
\IR_{f} M:=\sup_{x\in M} \rho_{\partial M,f}(x).
\end{equation*}
Having Proposition \ref{prop:global radial p-Laplacian comparison} at hand,
we can show the following comparison geometric result by the same argument as in the proof of \cite[Theorem 1.4 and Lemma 8.2]{Sak} for $N=1$,
whose unweighted version has been established by Kasue \cite{K3} for $p=2$, and the second named author \cite{Sa2} for general $p$ (cf. \cite[Theorem 2.1]{K3} and \cite[Theorem 1.6]{Sa2}):
\begin{thm}\label{thm:eigenvalue rigidity}
Let $p\in ]1,+\infty[$.
Assume
\begin{align*}
\Ric_f^N\geq c^{-1}\,\kappa\, e^{-\frac{4(1-\varepsilon)f}{n-1}}g,\quad H_{f,\partial M}\geq c^{-1}\lambda e^{-\frac{2(1-\varepsilon)f}{n-1}},\quad (1-\eps)f\leq (n-1)\delta
\end{align*}
for $\delta \in \mathbb{R}$.
Let $M$ be compact,
and let $f$ be $\partial M$-radial.
For $D\in ]0,C_{\kappa,\lambda}]\setminus \{+\infty\}$,
suppose $\IR_{f} M \leq D$.
Then
\begin{equation*}\label{eq:eigenvalue rigidity}
\nu_{\left( 1+\frac{2(1-\eps)}{n-1}  \right)f,p}(M)\geq \nu_{p,\kappa\,e^{-4\delta},\lambda\,e^{-2\delta},D\,e^{2\delta}}.
\end{equation*}
If the equality holds,
then $M$ is a $\left(\kappa e^{-4\delta} ,\lambda e^{-2\delta} \right)$-equational model space,
and $(1-\eps)f=(n-1)\delta$ on $M$.
\end{thm}

From the same method as in the proof of \cite[Lemma 8.4 and Theorem 8.8]{Sak} for $N=1$ with Proposition \ref{prop:global finite p-Laplacian comparison} and Theorem \ref{thm:finite inscribed radius rigidity},
we can also conclude the following:
\begin{thm}\label{thm:finite eigenvalue rigidity}
Let $p\in ]1,+\infty[$.
Let $\kappa$ and $\lambda$ satisfy the convex-ball-condition.
Assume
\begin{align*}
\Ric_f^N\geq c^{-1}\,\kappa\, e^{-\frac{4(1-\varepsilon)f}{n-1}}g,\quad H_{f,\partial M}\geq c^{-1}\lambda e^{-\frac{2(1-\varepsilon)f}{n-1}},\quad (1-\eps)f\leq (n-1)\delta
\end{align*}
for $\delta \in \mathbb{R}$.
Let $M$ be compact.
Then
\begin{equation*}
\nu_{f,p}(M)\geq \nu_{0,p}(B^{n}_{\kappa e^{-4\delta},\lambda e^{-2\delta}}).
\end{equation*}
If the equality holds,
then $(1-\eps)f=(n-1)\delta$ on $M$,
and the following properties hold:
\begin{enumerate}\setlength{\itemsep}{+0.7mm}
\item If $N=n$, then $M$ is isometric to $B^{n}_{\kappa\,e^{-4\delta},\lambda\,e^{-2\delta}}$;
\item if $N\neq n$, then $\eps =0$, and $M$ is isometric to $B^{n}_{\kappa\,e^{-4\delta},\lambda\,e^{-2\delta}}$.
\end{enumerate}
\end{thm}

\subsection{Spectrum rigidity}\label{sec:Spectrum rigidity}
Let $\Omega$ be a relatively compact domain in $M$ such that
its boundary is a smooth hypersurface in $M$ with $\partial \Omega \cap \partial M=\emptyset$.
For the canonical measure $v_{\partial \Omega}$ on $\partial \Omega$,
we set
\begin{equation*}
m_{f,\partial \Omega}:=e^{-f|_{\partial \Omega}}\, v_{\partial \Omega}.
\end{equation*}
Proposition \ref{prop:global finite p-Laplacian comparison} together with a similar calculation to the proof of \cite[Lemma 8.9]{Sak} for $N=1$ yields the following volume estimate,
which has been shown by Kasue \cite{K4} in the unweighted case (cf. \cite[Proposition 4.1]{K4}):
\begin{lem}\label{lem:Kasue volume estimate}
Let $\kappa$ and $\lambda$ satisfy the monotone-condition.
Assume
\begin{align*}
\Ric_f^N\geq c^{-1}\,\kappa\, e^{-\frac{4(1-\varepsilon)f}{n-1}}g,\quad H_{f,\partial M}\geq c^{-1}\lambda e^{-\frac{2(1-\varepsilon)f}{n-1}},\quad (1-\eps)f\leq (n-1)\delta
\end{align*}
for $\delta \in \mathbb{R}$.
Define $\rho_{\partial M,\delta}:=e^{-2\delta}\,\rho_{\partial M}$.
Let $\Omega$ be a relatively compact domain in $M$ such that
$\partial \Omega$ is a smooth hypersurface in $M$ satisfying $\partial \Omega \cap \partial M=\emptyset$.
Set
\begin{equation*}\label{eq:diameter of Omega}
D_{\delta,1}(\Omega):=\inf_{x\in \Omega}\, \rho_{\partial M,\delta}(x),\quad D_{\delta,2}(\Omega):=\sup_{x\in \Omega} \,\rho_{\partial M,\delta}(x).
\end{equation*}
Then we have
\begin{equation*}\label{eq:Kasue volume estimate}
m_{f}( \Omega) \leq e^{2\delta}\,  \sup_{s\in ]D_{\delta,1}(\Omega),D_{\delta,2}(\Omega)[}\,  \frac{\int^{ D_{\delta,2}(\Omega)}_{s}\,  \mathfrak{s}^{c^{-1}}_{\kappa,\lambda}(\xi)\, \d \xi}{\mathfrak{s}^{c^{-1}}_{\kappa,\lambda}(s)} \, m_{f,\partial \Omega}\, (\partial \Omega).
\end{equation*}
\end{lem}

For $D\in ]0,C_{\kappa,\lambda}]$,
we set
\begin{equation*}\label{eq:spectrum constant}
C(\kappa,\lambda,D):=\sup_{s\in [0,D[}\, \frac{\int^{D}_{s}\,  \mathfrak{s}^{c^{-1}}_{\kappa,\lambda}(\xi)\, \d \xi}{\mathfrak{s}^{c^{-1}}_{\kappa,\lambda}(s)}.
\end{equation*}
Notice that 
$C(\kappa,\lambda,+\infty)$ is finite if and only if $\kappa<0$ and $\lambda=\sqrt{\vert \kappa \vert}$;
in this case,
\begin{equation*}
C(\kappa,\lambda,D)=\left(c^{-1}\lambda \right)^{-1}\,\bigl(1-e^{-c^{-1}\lambda\, D} \bigl).
\end{equation*}
Lemma \ref{lem:Kasue volume estimate} and the argument in the proof of \cite[Lemma 8.11]{Sak} tell us the following:
\begin{lem}\label{lem:p-Laplacian1}
Let $p\in ]1,+\infty[$.
Let $\kappa$ and $\lambda$ satisfy the monotone-condition.
Assume
\begin{align*}
\Ric_f^N\geq c^{-1}\,\kappa\, e^{-\frac{4(1-\varepsilon)f}{n-1}}g,\quad H_{f,\partial M}\geq c^{-1}\lambda e^{-\frac{2(1-\varepsilon)f}{n-1}},\quad (1-\eps)f\leq (n-1)\delta
\end{align*}
for $\delta \in \mathbb{R}$.
For $D\in]0,C_{\kappa,\lambda}]$,
suppose $\IR M \leq e^{2\,\delta}\,D$.
Then we have
\begin{equation*}
\nu_{f,p}(M)\geq (\,p\,e^{2\delta}\,C(\kappa,\lambda,D)\,)^{-p}.
\end{equation*}
\end{lem}

Due to Theorem \ref{thm:splitting theorem},
Lemma \ref{lem:p-Laplacian1} and the argument in the proof of \cite[Theorem 8.12]{Sak} for $N=1$ lead us to the following:
\begin{thm}\label{thm:spectrum rigidity}
Let $p\in ]1,+\infty[$.
Let $\kappa<0$ and $\lambda:=\sqrt{\vert \kappa \vert}$.
Assume
\begin{align*}
\Ric_f^N\geq c^{-1}\,\kappa\, e^{-\frac{4(1-\varepsilon)f}{n-1}}g,\quad H_{f,\partial M}\geq c^{-1}\lambda e^{-\frac{2(1-\varepsilon)f}{n-1}},\quad (1-\eps)f\leq (n-1)\delta
\end{align*}
for $\delta \in \mathbb{R}$.
Let $\partial M$ be compact.
Then
\begin{equation*}\label{eq:spectrum rigidity}
\nu_{f,p}(M)\geq e^{-2p\delta}\left(\frac{c^{-1}\lambda }{p}\right)^{p}.
\end{equation*}
If the equality holds,
then $M$ is diffeomorphic to $[0,+\infty[\times \partial M$,
and the following hold:
\begin{enumerate}\setlength{\itemsep}{+0.7mm}
\item If $N=n$, then $f$ is constant, and
\begin{equation*}
g=dt^2+\s^2_{\kappa\,e^{-\frac{4(1-\eps)f}{n-1}},\lambda\,e^{-\frac{2(1-\eps)f}{n-1}}}(t)g_{\partial M};
\end{equation*} \label{enum:twisted curv}
\item if $N\neq 1,n$, then for any $z\in \partial M$
\begin{equation*}
f(\gamma_z(t))=f(z)-\eps \frac{N-n}{N-1}c^{-1}\log \mathfrak{s}_{\kappa,\lambda}(s_{f,z}(t)),\quad g=dt^2+\s^{2\frac{c^{-1}}{n-1} \left( 1-\eps \frac{N-n}{N-1}   \right) }_{\kappa,\lambda}(s_{f,z}(t))g_{\partial M};
\end{equation*}\label{enum:relaxed twisted curv}
\item if $N=1$, then
\begin{equation*}
\eps=0,\quad g=dt^2+\exp\left(2 \frac{f(\gamma_{z}(t))-f(z)}{n-1} \right)\,\s^2_{\kappa,\lambda}(s_{f,z}(t))\,g_{\partial M}.
\end{equation*} \label{enum:curv cond}
\end{enumerate}
\end{thm}

\noindent
\emph{Acknowledgment.}
The authors would thank Professor Xiang-Dong Li for his comment on the references.

\providecommand{\bysame}{\leavevmode\hbox to3em{\hrulefill}\thinspace}
\providecommand{\MR}{\relax\ifhmode\unskip\space\fi MR }
\providecommand{\MRhref}[2]{%
  \href{http://www.ams.org/mathscinet-getitem?mr=#1}{#2}
}
\providecommand{\href}[2]{#2}

\end{document}